\documentclass[11pt,letterpaper]{amsart}
\usepackage{amsfonts, amsmath, amssymb, amscd, amsthm}

\newtheorem{thm}{Theorem}

\newtheorem{lem}[thm]{Lemma}
\newtheorem{prop}[thm]{Proposition}

\theoremstyle{definition}
\newtheorem{defn}[thm]{Definition}
\theoremstyle{remark}

\renewcommand{\phi }{{\rm\bf Lab\, }}
\newcommand{\Ga }{\Gamma (G, \mathcal A)}
\newcommand{\G }{\Gamma (G, X\cup\mathcal H)}

\newcommand{\dxh }{dist_{X\cup\mathcal H}}
\newcommand{\dx }{dist_X}

\newcommand{\Hl }{\{ H_\lambda \} _{\lambda \in \Lambda }}
\newcommand{\e }{\varepsilon }
\newcommand{\ad }{{\rm asdim\, }}
\newfont{\eufm}{eufm10}

\begin{document}

\title{Asymptotic dimension of relatively hyperbolic groups}

\author{D. Osin}
\thanks{This work has been supported by the RFBR Grants
$\sharp $ 02-01-00892.}

\date{}



\address{Department of Mathematics, 1326 Stevenson Center,
Vanderbilt University, Nashville  TN 37240-0001, USA}

\email{denis.osin@gmail.com }

\begin{abstract}
Suppose that a finitely generated group $G$ is hyperbolic relative
to a collection of subgroups $\{ H_1, \ldots , H_m\} $. We prove
that if each of the subgroups $H_1, \ldots , H_m$ has finite
asymptotic dimension, then asymptotic dimension of $G$ is also
finite.
\end{abstract}

\maketitle


\section{Introduction}


The notion of asymptotic dimension of a metric space was proposed
by Gromov \cite{Gro} as a large--scale analog of the Lebesgue
covering dimension. Recall that a metric space $S$ has {\it
asymptotic dimension} $\ad S\le n$ if for any $r>0$, there exists
a covering $$S=\bigcup\limits_{\alpha\in A} U_\alpha
$$ such that the sets $U_\alpha $ are uniformly bounded and no more
than $n+1$ elements of $\{U_\alpha  \} _{\alpha \in A} $ meet any
ball of radius $r$.

In the case of finitely generated groups endowed with word
metrics, the question of finiteness of asymptotic dimension took
on additional significance with a theorem of Yu stating that the
Novikov Higher Signature Conjecture is true for manifolds $M$ such
that $\ad \pi _1(M)<\infty $ \cite{Yu}. Some other results
concerning groups of finite asymptotic dimension can be found in
\cite{CG,D,Yu1}. Although, in general, $\ad G$ can be infinite for
a finitely generated or even finitely presented group $G$, there
are many classes of groups for which asymptotic dimension is known
to be finite. For example, this is so for nilpotent groups,
fundamental groups of finite graphs of groups where vertex groups
have finite asymptotic dimension, hyperbolic groups, etc. (see
\cite{BD,BD2,R} and references therein).

In the present paper we study the case of relatively hyperbolic
groups. Originally the notion of relative hyperbolicity was
suggested by Gromov in \cite{Gro}. Since then it has been
elaborated from various points of view \cite{Bow,SD,F,RHG,Yaman}.
We recall the definitions of relative hyperbolicity suggested in
\cite{RHG}. In the case of finitely generated groups this
definition is equivalent to the definitions given in
\cite{Bow,F,SD,Yaman}.

Let $G$ be a group, $\Hl $ a collection of subgroups of $G$, $X$ a
subset of $G$. We say that $X$ is a {\it relative generating set
of $G$ with respect to $\Hl $} if $G$ is generated by $X$ together
with the union of all $H_\lambda $. (For convenience, we always
assume $X=X^{-1}$.) In this situation the group $G$ can be
regarded as a quotient of the free product
\begin{equation}
F=\left( \ast _{\lambda\in \Lambda } H_\lambda  \right) \ast F(X),
\label{F}
\end{equation}
where $F(X)$ is the free group with the basis $X$. Let $\e $
denote the natural homomorphism $F\to G$. If $Ker \, \e $ is a
normal closure of a subset $\mathcal R\subseteq N$ in the group
$F$, we say that $G$ has {\it relative presentation}
\begin{equation}\label{G}
\langle X,\; H_\lambda, \lambda\in \Lambda \; |\; R=1,\,
R\in\mathcal R \rangle .
\end{equation}
If $\sharp\, X<\infty $ and $\sharp\, \mathcal R<\infty $, the
relative presentation (\ref{G}) is said to be {\it finite} and the
group $G$ is said to be {\it finitely presented relative to the
collection of subgroups $\Hl $.}

Set
\begin{equation}\label{H}
\mathcal H=\bigsqcup\limits_{\lambda\in \Lambda} (H_\lambda
\setminus \{ 1\} ) .
\end{equation}
Given a word $W$ in the alphabet $X\cup \mathcal H$ such that $W$
represents $1$ in $G$, there exists an expression
\begin{equation}
W=_F\prod\limits_{i=1}^k f_i^{-1}R_i^{\pm 1}f_i \label{prod}
\end{equation}
with the equality in the group $F$, where $R_i\in \mathcal R$ and
$f_i\in F $ for $i=1, \ldots , k$. The smallest possible number
$k$ in a representation of the form (\ref{prod}) is called the
{\it relative area} of $W$ and is denoted by $Area^{rel}(W)$.

\begin{defn}
A group $G$ is {\it hyperbolic relative to a collection of
subgroups} $\Hl $ if $G$ is finitely presented relative to $\Hl $
and there is a constant $L>0$ such that for any word $W$ in $X\cup
\mathcal H$ representing the identity in $G$, we have $Area^{rel}
(W)\le L\| W\| $.
\end{defn}

We note that the above definition does not require the group $G$
and the subgroups $H_\lambda $ to be finitely generated as well as
the collection $\Hl $ to be finite. However, in case $G$ is
generated by a finite set in the ordinary (non--relative) sense
and is finitely presented relative to a collection of subgroups
$\Hl $, $\Lambda $ is known to be finite and the subgroups
$H_\lambda $ are known to be finitely generated \cite[Theorem
1.1]{RHG}.

The class of relatively hyperbolic groups includes many examples
of interest such as fundamental groups of finite-volume Riemannian
manifolds of pinched negative curvature \cite{Bow,F},
geometrically finite convergence groups \cite{Yaman}, small
cancellation quotients of free products \cite{RHG}, fully
residually free groups \cite{Dah}, etc. The main result of our
paper is the following.

\begin{thm}\label{main}
Suppose that a finitely generated group $G$ is hyperbolic relative
to a (finite) collection of subgroups $\Hl $ and each of the
subgroups $H_\lambda $ has finite asymptotic dimension. Then
asymptotic dimension of $G$ is finite.
\end{thm}

In the particular case when $G$ is hyperbolic relative to a
collection of virtually nilpotent subgroups $\Hl $, Theorem
\ref{main} was independently proved by Dahmani and Yaman
\cite{DY}. However their method essentially uses the assumption
about $H_\lambda $'s and can not be applied in the general case.

Recall that a group $G$ is said to be {\it weekly hyperbolic
relative to a collection of subgroups $\Hl $}, if the Cayley graph
$\G $ of $G$ with respect to the generating set $X\cup \mathcal H$
is hyperbolic, where $X$ is a finite generating set of $G$ modulo
$\Hl $ and $\mathcal H$ is the set defined by (\ref{H}). It is not
hard to show that relative hyperbolicity implies weak relative
hyperbolicity with respect to the same collection of subgroups.
The converse is not true. For instance, $G=\mathbb Z\times \mathbb
Z$ is weakly hyperbolic but not hyperbolic relative to the
multiples.

The main idea of the proof of Theorem \ref{main} is based on
exploring the weak relative hyperbolicity of $G$. On the other
hand, we also use certain additional arguments that do not follow
from the hyperbolicity of $\G $. So the natural question is
whether these arguments are essential.

In many particular cases, the weak hyperbolicity of $G $ relative
to a finite collection of subgroups $\Hl $ ensures the finiteness
of $\ad G$ whenever $\ad H_\lambda $ is finite for all $\lambda
\in \Lambda $. For example, if $H$ is a normal subgroup of $G$ and
$G/H$ is an ordinary hyperbolic group, then $G$ is weakly
hyperbolic relative to $H$ as $\Gamma (G, X\cup H)$ is
quasi--isometric to the quotient group $G/H$. In these settings,
$\ad H<\infty $ implies $\ad G<\infty$ by the result of Bell and
Dranishnikov \cite{BD2} stating that any extension of a group of
finite asymptotic dimension by a group of finite asymptotic
dimension has finite asymptotic dimension.

Another series of examples is provided by amalgamated free
products and HNN--extensions (or, more generally, by fundamental
groups of finite graphs of groups). Indeed any group of the form
$G=H_1\ast _{A=B}H_2$ is weakly hyperbolic relative to the
collection $\{ H_1, H_2\} $ \cite{WH}. In this case the finiteness
of asymptotic dimensions of $H_1$ and $H_2$ implies $\ad G<\infty
$ according to the main result of \cite{BD}. Similarly any
HNN--extension $G=H\ast _A$ is weakly hyperbolic relative to $H$
\cite{WH}. And again $\ad G<\infty $ whenever $\ad H<\infty $
\cite{BD}.

Note that, in general, the above--mentioned weakly relatively
hyperbolic groups are not relatively hyperbolic with respect to
the specified collections of subgroups. Generalizing these
examples one may conjecture that if a group $G$ is weakly
hyperbolic relative to a finite collection of subgroups $\Hl $ and
all subgroups $H_\lambda $ have finite asymptotic dimension, then
asymptotic dimension of $G$ is finite. However this conjecture
does not hold. To provide a counterexample we prove the following.

\begin{prop}\label{prop}
There exists a finitely presented group $G$ and a finite
collection of subgroups $\Hl $ such that:

\begin{enumerate}
\item[1)] For any $\lambda \in \Lambda $, $H_\lambda $ is cyclic
and hence $\ad H_\lambda \le 1$.

\item[2)] The Cayley graph $\G $ has finite diameter; in
particular, it is hyperbolic and thus $G$ is weakly hyperbolic
relative to $\Hl $.

\item[3)] $\ad G=\infty $.
\end{enumerate}
\end{prop}

The paper is organized as follows. In the next section we collect
all necessary definitions and results used in our paper. Sections
3 and 4 contain proofs of certain technical facts about relatively
hyperbolic groups. The proofs of Theorem \ref{main} and
Proposition \ref{prop} are provided in Section 5 and 6
respectively.


\section{Preliminaries}


\subsection{Some notation and conventions}
Given a word $W$ in an alphabet $\mathcal A$, we denote by $\| W\|
$ its length that is the number of letters in $W$. We also write
$W\equiv V$ to express the letter for letter equality of words $W$
and $V$. For elements $g$, $t$ of a group $G$, $g^t$ denotes the
element $t^{-1}gt$. Similarly $H^t$ denotes $t^{-1}Ht$ for a
subgroups $H\le G$.  Recall that a subset $X$ of a group $G$ is
said to be {\it symmetric} if for any $x\in X$, we have $x^{-1}\in
X$. In this paper all generating sets of groups under
consideration are supposed to be symmetric.

Given a group $G$ generated by a (symmetric) set $\mathcal A$, the
{\it Cayley graph} $\Ga $ of $G$ with respect to $\mathcal A$ is
an oriented labelled 1--complex with the vertex set $V(\Ga )=G$
and the edge set $E(\Ga )=G\times \mathcal A$. An edge $e=(g,a)$
goes from the vertex $g$ to the vertex $ga$ and has label $\phi
(e)\equiv a$. As usual, we denote the origin and the terminus of
the edge $e$, i.e., the vertices $g$ and $ga$, by $e_-$ and $e_+$
respectively. Given a combinatorial path $p=e_1e_2\ldots e_k$ in
the Cayley graph $\Ga $, where $e_1, e_2, \ldots , e_k\in E(\Ga
)$, we denote by $\phi (p)$ its label. By definition, $$\phi
(p)\equiv \phi (e_1)\phi(e_2)\cdots \phi (e_k).$$ We also denote
by $p_-=(e_1)_-$ and $p_+=(e_k)_+$ the origin and the terminus of
$p$ respectively. The length $l(p)$ of $p$ is the number of edges
of $p$.

Associated to $\mathcal A$ is the so--called {\it word metric} on
$G$. More precisely, the length $|g|_\mathcal A$ of an element
$g\in G$ is defined to be the length of a shortest word in
$\mathcal A$ representing $g$ in $G$. This defines a metric on $G$
by $dist_\mathcal A(f,g)=|f^{-1}g|_\mathcal A$. We also denote by
$dist _\mathcal A$ the natural extension of the word metric to the
Cayley graph $\Ga $.

\subsection{Relatively hyperbolic groups}

Recall that a metric space $M$ is {\it $\delta $--hyperbolic} for
some $\delta \ge 0$ (or simply {\it hyperbolic}) if for any
geodesic triangle in $M$, any side of the triangle belongs to the
union of the closed $\delta $--neighborhoods of the other two
sides. A group $G$ is called (ordinary) {\it hyperbolic} if $G$
satisfies Definition 1 with respect to the trivial subgroup. An
equivalent definition says that $G$ is hyperbolic if it is
generated by a finite set $X$ and the Cayley graph $\Gamma (G, X)$
is a hyperbolic metric space. In the relative case these
approaches are not equivalent, but we still have the following
\cite[Theorem 1.7]{RHG}.

\begin{lem}\label{CG}
Suppose that $G$ is a group hyperbolic relative to a collection of
subgroups $\Hl $. Let $X$ be a finite relative generating set of
$G$ with respect to $\Hl $. Then the Cayley graph $\G $ of $G$
with respect to the generating set $X\cup \mathcal H$ is a
hyperbolic metric space.
\end{lem}

Let us recall an auxiliary terminology introduced in \cite{RHG},
which plays an important role in our paper. Let $G$ be a group,
$\Hl $ a collection of subgroups of $G$, $X$ a finite generating
set of $G$ with respect to $\Hl $, $q$ a path in the Cayley graph
$\G $. A subpath $p$ of $q$ is called an {\it $H_\lambda
$--component} for some $\lambda \in \Lambda $ (or simply a {\it
component}) of $q$, if the label of $p$ is a word in the alphabet
$H_\lambda\setminus \{ 1\} $ and $p$ is not contained in a bigger
subpath of $q$ with this property.

Two components $p_1, p_2$ of a path $q$ in $\G $ are called {\it
connected} if they are $H_\lambda $--components for the same
$\lambda \in \Lambda $ and there exists a path $c$ in $\G $
connecting a vertex of $p_1$ to a vertex of $p_2$ such that ${\phi
(c)}$ entirely consists of letters from $ H_\lambda $. In
algebraic terms this means that all vertices of $p_1$ and $p_2$
belong to the same coset $gH_\lambda $ for a certain $g\in G$.
Note that we can always assume $c$ to have length at most $1$, as
every nontrivial element of $H_\lambda $ is included in the set of
generators.  An $H_\lambda $--component $p$ of a path $q$ is
called {\it isolated } if no distinct $H_\lambda $--component of
$q$ is connected to $p$.

If the group $G$ is generated by the finite set $X$ in the
ordinary (non--relative) sense, we have two word metrics $\dx $
and $\dxh $ on $G$ associated to the generating sets $X$ and
$X\cup \mathcal H$ respectively. Each of these has its own
advantage and disadvantage. Namely, if $G$ is hyperbolic relative
to $\Hl $, then the metric space $(G, \dxh )$ is hyperbolic, but
in general is not locally finite. On the other hand, $(G, \dx )$
is locally finite but usually is not hyperbolic. In the next two
lemmas we consider these metric spaces together. The first result
is a simplification of Lemma 2.27 from \cite{RHG}.

\begin{lem}\label{Omega}
Suppose that a group $G$ is generated by a finite set $X$ and is
hyperbolic relative to a collection of subgroups $\Hl $. Then
there exists a constant $L>0$ such that for any cycle $q$ in $\G
$, any $\lambda\in \Lambda $, and any set of isolated $H_\lambda
$--components $p_1, \ldots , p_k$ of $q$, we have
$$ \sum\limits_{i=1}^k \dx ((p_i)_-, (p_i)_+)\le Ll(q).$$
\end{lem}

Note that if $p$ is a geodesic path in $\G $, then each component
of $p$ is isolated and consists of a single edge. The following
lemma is a particular case of Theorem 3.23 from \cite{RHG}. (In
Farb's approach \cite{F}, this is a part of the definition of a
relatively hyperbolic group called the Bounded Coset Penetration
property.)

\begin{lem}\label{geod}
Suppose that a group $G$ is generated by a finite set $X$ and is
hyperbolic relative to a collection of subgroups $\Hl $. Then for
any $s\ge 0$, there exists a constant $\e =\e (s) \ge 0$ such that
the following condition holds. Let $p_1, p_2$ be two geodesics in
$\G $ such that $$\max \{ \dx ((p_1)_-, (p_2)_-),\, \dx ((p_1)_+,
(p_2)_+)\} \le s$$ and let $c$ be a component of $p_1$ such that
$\dx (c_-, c_+)\ge \e $. Then there is a component of $p_2$
connected to $c$.
\end{lem}

\subsection{Asymptotic dimension}

We also recall some properties of asymptotic dimension used in our
paper.  The first property is quite obvious and immediately
follows from the definition.

\begin{lem}\label{ad1}
If $M_1$ is a metric space and $M_2\subseteq M_1$ is a subspace
endowed with the induced metric, then $\ad  M_2 \le \ad  M_1 $.
\end{lem}

Recall that a map $\alpha : M_1\to M_2$ from a metric space $M_1$
to a metric space $M_2$ is called a {\it quasi--isometry}, if
there are constants $\lambda >0$, $c\ge 0$, $\e \ge 0$ such that
the following conditions hold:

\begin{enumerate}
\item For any two points $x,y\in M_1$, we have $$ \frac1\lambda
dist_{M_2} (\alpha (x), \alpha (y))-c\le dist_{M_1} (x, y)\le
\lambda dist_{M_2} (\alpha (x), \alpha (y))+c .$$

\item The image $\alpha (M_1)$ is $\e $--dense in $M_2$, that is,
for any $z\in M_2$ there exists $x\in M_1$ such that $dist_{M_2}
(\alpha (x), z)\le \e $.
\end{enumerate}

Two metric spaces $M_1 $ and $M_2$ are called {\it
quasi--isometric}, if there exists a quasi--isometry from $M_1 $
to $M_2$. It is not hard to check that this is an equivalence
relation. The lemma below is also quite obvious (see, for example,
\cite[Sec. 9.1]{R}).

\begin{lem}\label{qi}
If $M_1$ and $M_2$ are quasi--isometric, then $\ad M_1=\ad M_2$.
\end{lem}

The next two results were proved by Bell and Dranishnikov
\cite{BD1}.

\begin{lem} \label{ad2}
Let $M=M^\prime \cup M^{\prime\prime } $ be a metric space. Then
$$\ad M =\max \{ \ad M^\prime , \, \ad M^{\prime \prime }\} .$$
\end{lem}

As pointed out in \cite{BD1}, Lemma \ref{ad2} can be generalized
to certain infinite unions. More precisely, one says that a
collection of spaces $\{ M_\alpha \} _{\alpha \in A} $ has {\it
asymptotic dimension $\le n$ uniformly}, if for any $r>0$, there
exist coverings $\mathcal U_\alpha $ of $M_\alpha $ and a constant
$d$ such that for all $\alpha $, all $U\in \mathcal U_\alpha $
have diameter at most $d$ and any ball of radius $r$ in $M_\alpha
$ intersects at most $n+1$ elements of the covering $\mathcal
U_\alpha $. Recall also that two subsets $A$, $B$ in a metric
space are called {\it $s$--separated} for some $s>0$ if
$dist(a,b)\ge s$ for all $a\in A$, $b\in B$.

\begin{lem} \label{ad3}
Suppose that $M$ is a metric space and $M=\bigcup\limits_{\alpha }
M_\alpha $, where $M_\alpha $ have asymptotic dimension at most
$n$ uniformly. Suppose also that for any $s>0$, there is $Y_s\in
M$ such that $\ad Y_s\le n$ and the sets $M_\alpha \setminus Y_s$
(fixed $s$, varying $\alpha $) are $s$--separated. Then $\ad M\le
n$.
\end{lem}

Finally, let $G$ be a group acting on a matric space $M$. Given
$x\in M$ and $R\ge 0$, we define an $R$--quasi--stabilizer of $x$
by $$W_R(x)=\{ g\in G\; |\; dist (x, gx)\le R\} .$$ The lemma
below appears as Theorem 2 in \cite{BD1}.

\begin{lem} \label{ad4}
Suppose that a finitely generated group $G$ acts by isometries on
a metric space $M$ such that $\ad M<m $ and for some $x\in M$,
$\ad W_R(x)\le n $. Then $\ad G\le (m+1)(n+1)-1 $.
\end{lem}


\section{Asymptotic dimension of relative balls}


Throughout the rest of the paper, $G$ denotes a group that is
generated by a finite set $X$ and is hyperbolic relative to a
collection of subgroups $\Hl $, where $\sharp\, \Lambda <\infty $.
We keep the notation for Cayley graphs, word metrics, etc.,
introduced in the previous section. Speaking on asymptotic
dimension of subsets of the group $G$ we always mean the
asymptotic dimension associated to the metric $\dx $.

We begin with an auxiliary result.

\begin{lem}\label{lc}
Let $p_1=q_1e_1$, $p_2=q_2e_2$ be two geodesics in $\G $ such that
$$\max \{ \dx ((p_1)_-, (p_2)_-),\, \dx ((p_1)_+, (p_2)_+)\} \le
s.$$ Suppose that for a certain $\lambda\in \Lambda $,  $e_1$ and
$e_2$ are $H_\lambda $--components satisfying the inequality
\begin{equation}\label{dei}
\dx ((e_i)_-, (e_i)_+)\ge \max \{ \e ,\, 2L(s+1)\} ,
\end{equation}
where $\e =\e (s)$ and $L$ are constants provided by Lemma
\ref{geod} and Lemma \ref{Omega} respectively. Then $e_1$ and
$e_2$ are connected.
\end{lem}

\begin{figure}
 \unitlength 1mm 
\linethickness{0.4pt}
\ifx\plotpoint\undefined\newsavebox{\plotpoint}\fi 
\begin{picture}(97.54,41)(5,9)
\put(44.01,38.66){\vector(4,-1){.07}}\qbezier(6.01,47.5)(47.55,35.84)(74.95,35.48)
\qbezier(6.09,26.76)(40.58,16.43)(75.07,14.72)

\linethickness{1pt}
\multiput(75.07,35.58)(4.4,-.03){5}{\line(1,0){4.4}}
\qbezier(33,19.88)(44.25,17.94)(55,16.25)
\put(75.13,14.63){\line(1,0){22}}
\qbezier(55.13,16.38)(69.69,26.56)(97,35.5)

\thicklines

\put(55.07,12.6){$t_+$}\put(98.8,14.55){$(p_2)_+$}
\put(98.3,35.45){$(p_1)_+$}

\put(86.07,35.5){\vector(1,0){.07}} \put(75.07,35.58){\circle*{1}}
\put(96.94,35.45){\circle*{1}} \put(6.07,47.58){\circle*{1}}
\put(75.07,14.55){\circle*{1}} \put(97.07,14.55){\circle*{1}}
\put(55.07,16.42){\circle*{1}} \put(33.07,19.92){\circle*{1}}
\put(6.07,26.8){\circle*{1}} \put(44.13,18){\vector(4,-1){.07}}
\put(86.13,14.63){\vector(1,0){.07}}
\put(42.38,41.63){\makebox(0,0)[cc]{$q_1$}}
\put(25.5,17.5){\makebox(0,0)[cc]{$q_2$}}
\put(44,20.8){\makebox(0,0)[cc]{$t$}}
\put(85.13,38){\makebox(0,0)[cc]{$e_1$}}
\put(85.25,11.75){\makebox(0,0)[cc]{$e_2$}} \thinlines
\put(95.1,23.5){\vector(0,-1){.07}}\qbezier(96.88,35.38)(93,22)(97.13,14.63)
\put(97.13,26){\makebox(0,0)[cc]{$u$}} \thicklines
\put(72.88,26.25){\vector(2,1){.07}} \put(70,27){$c$}
\end{picture}
  \caption{}\label{lc-fig}
\end{figure}

\begin{proof}
By Lemma \ref{geod}, there is a component $t$ of $q_2$ that is
connected to $e_1$. Assume that $t\ne e_2$. Let $c$ denote a path
in $\G $ of length at most $1$ labelled by an element of
$H_\lambda $ such that $e_-=t_+$, $e_+=(e_1)_+$ (see Fig.
\ref{lc-fig}). We also denote by $u$ a path connecting $(p_1)_+$
to $(p_2)_+$ such that $\phi (u)$ is a word in the alphabet $X$
and $l(u)\le s$. Consider the cycle $q=cu[(p_2)_+, t_+]$, where
$[(p_2)_+, t_+]$ is a segment of $p_2^{-1}$. Note that $e_2$ is an
isolated component of $q$. Indeed $e_2$ is isolated in $[(p_2)_+,
t_+]$ since $p_2$ is geodesic. Further $e_2$ is not connected to
$c$. Indeed otherwise $e_2$ is connected to $t$ that contradicts
the assumption that $p_2$ is geodesic again. Finally $e_2$ can not
be connected to a component of $u$ as $\phi (u)$ is a word in $X$
and thus $u$ contains no components at all.

Since $p_2$ is geodesic, we have $$ l([(p_2)_+, t_+]) \le
l(c)+l(u)\le s+1.$$ Consequently, $$ l(q)= l([(p_2)_+, t_+])
+l(c)+l(u) \le 2(s+1).$$ Applying Lemma \ref{Omega} we obtain the
inequality
$$\dx ((e_2)_-, (e_2)_+)\le Ll(q)\le 2L(s+1),$$ which contradicts
(\ref{dei}). The lemma is proved.
\end{proof}

Denote by $B(n)$ the ball in $G$ centered at $1$ of radius $n$
with respect to the metric $\dxh $, i.e., $$ B(n)=\{ g\in G\; | \;
|g|_{X\cup \mathcal H}\le n\} .$$ We have already pointed out that
$B(n)$ is not necessarily finite. The main result of this section
is the following.

\begin{lem}\label{Bn}
Suppose that all subgroups $H_\lambda $, $\lambda \in \Lambda $,
have asymptotic dimension at most $d$. Then for any $n\in \mathbb
N$, we have $\ad B(n)\le d$.
\end{lem}

\begin{proof}
We proceed by induction on $n$. For $n=1$, we have $B(1)=X\cup
\left(\bigcup\limits_{\lambda \in \Lambda } H_\lambda \right) $.
Since $X$ and $\Lambda $ are finite, the inequality $\ad B(1)\le
d$ follows from Lemma \ref{ad2}.

Now assume $n>1$. Then clearly $$ B(n)= \left(
\bigcup\limits_{\lambda \in \Lambda } B(n-1) H_\lambda \right)
\cup \left( \bigcup\limits_{x\in X } B(n-1) x \right) .
$$
Note that the identity map $B(n-1)\to B(n-1)$ induces a
quasi--isometry from $B(n-1)$ to $B(n-1)x$. Hence, $$\ad \left(
\bigcup\limits_{x\in X } B(n-1) x \right) \le d$$ according to
Lemma \ref{qi}, Lemma \ref{ad2}, and the inductive assumption.
Thus it remains to show that
\begin{equation}\label{rest}
\ad B(n-1) H_\lambda  <d
\end{equation}
for any $\lambda \in \Lambda $.

Throughout the rest of the proof we fix an arbitrary $\lambda \in
\Lambda $ and denote by $R(n-1)$ the subset of $B(n-1)$ such that
for any $b\in B(n-1)$, we have $bH_\lambda =gH_\lambda $ for a
certain $g\in R(n-1)$ and $gH_\lambda \ne fH_\lambda $ whenever
$f$, $g$ are different elements of $R(n-1)$. Obviously we have $$
B(n-1)H_\lambda =\bigsqcup\limits_{g\in R(n-1)} gH_\lambda .$$

Let us fix some $s>0$ and set $$ T_s=\{ g\in G\; |\; |g|_X\le \max
\{ \e ,\, 2L(s+1)\} \} ,$$ where $\e =\e (s)$ and $L$ are the
constants from Lemma \ref{geod} and Lemma \ref{Omega}
respectively. Further we define $$Y_s=B(n-1)T_s.$$ Since $\sharp\,
T_s<\infty $, we have $\ad Y_s =\ad B(n-1)\le d$.  Let us show
that the sets $gH_\lambda \setminus Y_s$, $g\in R(n-1)$, are
$s$--separated.

Suppose that $x\in g_1H_\lambda \setminus Y_s$, $y\in g_2H_\lambda
\setminus Y_s$ for different $g_1, g_2\in R(n-1)$. Then
$x=g_1h_1$, $y=g_2h_2$ for some $h_1, h_2\in H_\lambda \setminus
T_s$. Assume that $\dx (x, y)\le s$. Let $A_i$, $i=1,2$, denote a
shortest word in $X\cup \mathcal H$ representing $g_i$ in $G$. Let
also $p_i$, $i=1,2$, denote the path in $\G $ such that
$(p_i)_-=1$ and $\phi (p_i)=A_ih_i$. Clearly $p_i$ is geodesic in
$\G $. Indeed otherwise we would have $$ |g_ih_i|_{X\cup \mathcal
H} =\dxh ((p_i)_-, (p_i)_+) < l(p_i)=\| A_i\| +1 =| g_i|_{X\cup
\mathcal H} +1=n$$ and hence $g_ih_i\in B(n-1)\subseteq Y_s$ that
contradicts our assumption. Note also that $$\dx ((p_1)_+,
(p_2)_+)=\dx (x,y)\le s.$$ As $h_i\notin T_s$, we have $|h_i|_X>
\max \{ \e ,\, 2L(s+1)\} $ for $i=1,2$. Therefore, the $H_\lambda
$--components of $p_1$ and $p_2$ labelled $h_1$ and $h_2$
respectively are connected by Lemma \ref{lc}. This means that
$g_1H_\lambda =g_2H_\lambda $ contradictory the choice of
$R(n-1)$.

Thus the sets $gH_\lambda \setminus Y_s$, $g\in R(n-1)$, are
$s$--separated. To complete the proof of (\ref{rest}) it remains
to apply Lemma \ref{ad3}.
\end{proof}


\section{Geodesic triangles in $\G $}


Let $\Delta =\Delta (x,y,z)$ be a geodesic triangle in a metric
space with vertices $x,y,z$. The triangle inequality tells us that
there exist (unique) points $a\in [y,z]$, $b\in [z,x]$, and $c\in
[x,y]$ such that $dist(x,b)=dist (x,c)$, $dist(y,a)=dist(y,c)$,
and $dist(z,a)=dist(z,b)$. Recall that two points $u\in [x,b]$ and
$v\in [x,c]$ are {\it conjugate} if $dist (x,u)=dist (x,v)$. In
the same way one defines conjugate points on pairs of segments
$[y, a]$, $[y, c]$ and $[z,a]$, $[z, b]$. The triangle $\Delta $
is said to be {\it $\xi $--thin} if $dist (u,v)\le \xi $ for any
two conjugate points $u,v\in \Delta $.

The following observation is quite obvious. We leave the proof to
the reader.

\begin{lem}\label{uv}
Let $\Delta (x,y,z)$ be a geodesic triangle in a metric space.
Suppose that $u\in [x,y]$, $v\in [x,z]$, $dist (x,u)=dist (x,v)$,
and $$ dist (u,y)+dist (v,z) \ge dist (y,z).$$ Then $u$ and $v$
are conjugate.
\end{lem}

The next lemma provides an equivalent definition of hyperbolicity.
(see, for example, \cite[Ch. III.H, Prop. 1.17]{BriH}).

\begin{lem}\label{thin}
A geodesic metric space $M$ is hyperbolic if and only if there
exists $\xi \ge 0$ such that every geodesic triangle in $M$ is
$\xi $--thin.
\end{lem}

In particular, any geodesic triangle in $\G $ is $\xi $--thin for
some $\xi =\xi (G)$. For our goals we need a stronger result
stating that geodesic triangles in $\G $ are thin with respect to
the metric $\dx $ associated to the finite generating set $X$. To
simplify our exposition we do not prove this result in the full
generality and restrict ourselves to the following particular
case, which is sufficient for our goal.

\begin{lem} \label{xi}
There exist constants $\rho ,\sigma >0$ having the following
property. Let $\Delta $ be a triangle with vertices $x$, $y$, $z$,
whose sides $[x,y]$, $[y,z]$, $[x,z]$ are geodesics in $\Gamma (G,
X\cup \mathcal H)$. Suppose that $u$ and $v$ are vertices on
$[x,y]$ and $[x,z]$ respectively such that
\begin{equation}\label{dxuxv}
\dxh (x,u)=\dxh (x,v)
\end{equation}
and
\begin{equation}\label{sigma}
\dxh (u,y)+\dxh (v, z) \ge \dxh (y,z) + \sigma .
\end{equation}
Then
\begin{equation}\label{rho}
\dx (u,v)\le \rho .
\end{equation}
\end{lem}

\begin{figure}
  \unitlength 1mm 
\linethickness{0.4pt}
\ifx\plotpoint\undefined\newsavebox{\plotpoint}\fi
\begin{picture}(102.13,50)(7,25)
\qbezier[1000](13.08,50.91)(73.19,59.57)(100.76,71.77)
\qbezier[1000](100.76,71.77)(87.86,44.9)(94.75,25.81)
\qbezier[1000](94.75,25.81)(72.66,43.75)(13.08,50.73)
\put(50.99,57.23){\line(0,-1){12.78}}
\qbezier(39.99,55)(42.81,49.5)(37.91,46.97)
\qbezier(63.18,59.91)(59.83,49.8)(63.92,40.88)
\put(50.99,53.81){\circle*{1}} \put(40.49,53.56){\circle*{1}}
\put(39.87,48.56){\circle*{1}} \put(50.86,47.19){\circle*{1}}
\put(50.86,44.44){\circle*{1}} \put(51.24,57.31){\circle*{1}}
\put(40.11,55.19){\circle*{1}} \put(37.99,47.19){\circle*{1}}
\put(63.99,40.94){\circle*{1}} \put(63.36,60.06){\circle*{1}}
\put(100.74,71.94){\circle*{1}} \put(94.74,25.81){\circle*{1}}
\put(12.99,50.69){\circle*{1}}

\thicklines \put(40.77,49.8){\vector(-1,-4){.07}}
\put(51,49){\vector(0,-1){.07}}

\put(12.25,48.13){\makebox(0,0)[cc]{$x$}}
\put(103,72.88){\makebox(0,0)[cc]{$y$}}
\put(96.25,24.25){\makebox(0,0)[cc]{$z$}}
\put(51,59.88){\makebox(0,0)[cc]{$u$}}
\put(50.88,41.63){\makebox(0,0)[cc]{$v$}}
\put(37.5,56.63){\makebox(0,0)[cc]{$u_1$}}
\put(35.75,45.25){\makebox(0,0)[cc]{$v_1$}}
\put(63,63){\makebox(0,0)[cc]{$u_2$}}
\put(64.25,36.75){\makebox(0,0)[cc]{$v_2$}}
\put(53.3,50.75){\makebox(0,0)[cc]{$s$}}

\thinlines \put(61.75,50.63){\vector(0,-1){.07}}
\put(61.63,51.13){\line(0,-1){.5}}

\linethickness{1pt} \put(51,53.88){\line(0,-1){6.63}}
\qbezier(40.38,53.5)(41.88,51.81)(39.88,48.38)
\put(38.88,51.5){\makebox(0,0)[cc]{$t$}}
\qbezier(40.5,53.63)(45.56,51.06)(50.88,53.75)

\end{picture}
  \caption{}\label{0}
\end{figure}

\begin{proof}
We set
\begin{equation}\label{s=4x}
\sigma =5\xi ,
\end{equation}
where $\xi $ is the constant provided by Lemma \ref{thin} for
$M=\G $. Note that $u$ and $v$ are conjugate according to Lemma
\ref{uv} and the inequality (\ref{sigma}). We denote by $p$ a
geodesic in $\Gamma $ such that $p_-=u$, $p_+=v$. By the choice of
$\xi $ we have
\begin{equation}\label{l1}
l(p)\le \xi .
\end{equation}

Further let $u_1\in [x,u]$, $v_1\in [x,v]$ be the vertices chosen
as follows. If $\dxh (x,u)=\dxh (x,v)< 2\xi $, we set $u_1=v_1=x$.
Otherwise $u_1$, $v_1$ are uniquely defined by the equality
\begin{equation}\label{l2}
\dxh (u_1, u)=\dxh (v_1, v) =2\xi .
\end{equation}
Obviously $u_1$ and $v_1$ are conjugate. Similarly let $u_2$ and
$v_2$ be the vertices on the segments $[u,y]$ and $[v,z]$
respectively such that
\begin{equation}\label{l3}
\dxh (u_2, u)=\dxh (v_2, v) =2\xi .
\end{equation}
Note that such vertices always exist. Indeed the inequalities
(\ref{sigma}) and (\ref{s=4x}) imply
$$
\begin{array}{rl}
\dxh (u,y)\ge & \frac12 \left( \dxh (u,y) +\dxh (v,z) -\dxh
(y,z)-\dxh(u,v)\right) \ge \\ & \\ & \frac12 (\sigma -\xi )=2 \xi
\end{array}
$$
and similarly $\dxh (v,z)\ge 2\xi $. The vertices $u_2$ and $v_2$
are conjugate by Lemma \ref{uv} since
$$ \begin{array}{rl}
\dxh (u_2,y)+\dxh (v_2, z)= & \dxh (u,y)+\dxh (v,z) -4 \xi \ge
\\ & \\ & \dxh (y,z) \end{array}$$
according to (\ref{sigma}) and (\ref{s=4x}). We denote by $o_1$
and $o_2$ geodesic paths in $\G $ such that $(o_1)_-=u_1$,
$(o_1)_+=v_1$, $(o_2)_-=u_2$, $(o_2)_+=v_2$. By Lemma \ref{xi} we
have
\begin{equation}\label{l4}
l(o_i)\le \xi  , \;\;\;\;\; i=1,2.
\end{equation}

Let us consider an arbitrary $H_\lambda $--component $s$ of $p$
for some $\lambda \in \Lambda $. In order to obtain an upper bound
on the distance $\dx (s_-, s_+)$, we are going to show that $s$ is
an isolated component in at least one of the cycles
$$
c_i=p[v,v_i]o_i^{-1}[u_i,u],\;\;\;\;\; i=1,2.
$$

First we note that $s$ can not be connected to a component of
$o_1$ or $o_2$. Indeed if, for example, $s$ is connected to a
component $t$ of $o_1$ (see Fig. \ref{0}), then $\dxh (t_-,
s_-)\le 1$. Taking into account (\ref{l4}) we obtain
$$
\begin{array}{rl}
\dxh (u_1, u)\le & \dxh (u_1, t_-) +\dxh (t_-, s_-) +\dxh (s_-,
u)\le \\ & \\ & (l(o_1)-1) +1 +(l(p)-1) < 2\xi .
\end{array}
$$
According to the choice of $u_1$ and $v_1$, this means
$u_1=v_1=x$. However, in this case $o_1$ is trivial and can not
contain components. A similar argument shows that $s$ can not be
connected to a component of $o_2$.

\begin{figure}
  \unitlength 1mm 
\linethickness{0.4pt}
\ifx\plotpoint\undefined\newsavebox{\plotpoint}\fi
\begin{picture}(102.13,50)(7,25)
\qbezier[1000](13.08,50.91)(73.19,59.57)(100.76,71.77)
\qbezier[1000](100.76,71.77)(87.86,44.9)(94.75,25.81)
\qbezier[1000](94.75,25.81)(72.66,43.75)(13.08,50.73)
\put(50.99,57.23){\line(0,-1){12.78}}
\qbezier(39.99,55)(42.81,49.5)(37.91,46.97)
\qbezier(63.18,59.91)(59.83,49.8)(63.92,40.88)
\put(50.99,53.81){\circle*{1}} \put(50.86,47.19){\circle*{1}}
\put(50.86,44.44){\circle*{1}} \put(51.24,57.31){\circle*{1}}
\put(40.11,55.19){\circle*{1}} \put(37.99,47.19){\circle*{1}}
\put(63.99,40.94){\circle*{1}} \put(63.36,60.06){\circle*{1}}
\put(100.74,71.94){\circle*{1}} \put(94.74,25.81){\circle*{1}}
\put(12.99,50.69){\circle*{1}} \put(42.49,55.81){\circle*{1}}
\put(48.11,56.81){\circle*{1}} \put(53.86,43.56){\circle*{1}}
\put(60.11,41.94){\circle*{1}} \linethickness{1pt}
\put(51,53.88){\line(0,-1){6.63}}
\multiput(42.38,55.75)(.191667,.033333){30}{\line(1,0){.191667}}
\multiput(53.75,43.5)(.1417778,-.0333333){45}{\line(1,0){.1417778}}
\put(12.25,48.13){\makebox(0,0)[cc]{$x$}}
\put(103,72.88){\makebox(0,0)[cc]{$y$}}
\put(96.25,24.25){\makebox(0,0)[cc]{$z$}}
\put(51,59.88){\makebox(0,0)[cc]{$u$}}
\put(50.88,41.63){\makebox(0,0)[cc]{$v$}}
\put(37.5,56.63){\makebox(0,0)[cc]{$u_1$}}
\put(35.75,45.25){\makebox(0,0)[cc]{$v_1$}}
\put(63,63){\makebox(0,0)[cc]{$u_2$}}
\put(64.25,36.75){\makebox(0,0)[cc]{$v_2$}}
\qbezier(42.38,55.75)(47.31,48.25)(51,47.25)
\put(53,51.5){\makebox(0,0)[cc]{$s$}}
\put(44.88,58.5){\makebox(0,0)[cc]{$a$}}
\put(56.38,40.4){\makebox(0,0)[cc]{$b$}}
\qbezier(50.88,47.25)(58.25,45.69)(60.13,41.88) \thinlines

\put(39.13,51){\makebox(0,0)[cc]{$o_1$}}
\put(64.5,51.75){\makebox(0,0)[cc]{$o_2$}}
\put(61.75,50.63){\vector(0,-1){.07}}
\put(41.1,50.63){\vector(0,-1){.07}}
 \thicklines
\put(51,49){\vector(0,-1){.07}}
\put(46.7,56.59){\vector(4,1){.07}}
\put(58.5,42.3){\vector(3,-1){.07}}
\end{picture}
  \caption{}\label{T}
\end{figure}

Further assume that there are edges $a$ and $b$ of the cycles
$c_1$ and $c_2$ respectively labelled by elements of $H_\lambda $
such that $a\ne s$, $b\ne s$, and $a_\pm $ and $b_\pm $ are
connected to $s_\pm $ by paths of lengths at most $1$ labelled by
elements of $H_\lambda $. As shown above, $a$, $b$ can not belong
to $o_1$ or $o_2$. Since $p$ is geodesic, $a$ and $b$ can not
belong to $p$. Similarly as $[x,y]$ and $[x,z]$ are geodesic, $a$
and $b$ can not belong to $[x,y]$ (or $[x,z]$) simultaneously.
Thus the only possibility is $a\in [x,y]$, $b\in [x,z]$ (or
conversely). For definiteness, we assume $a\in [x,y]$, $b\in
[x,z]$ (see Fig. \ref{T}). In this case we have
$$
\begin{array}{rl}
\dxh (x,v)\le & \dxh (x, b_+) -1 \le \dxh (x, a_-)+\dxh (a_-,
b_+)-1\le \\ &\\ & \dxh (x, a_-)< \dxh (x,u)
\end{array}
$$
contradictory (\ref{dxuxv}).

Thus each component of $p$ is isolated in at least one of the
cycles $c_1$, $c_2$. Combining (\ref{l1})--(\ref{l4}) yields $
l(c_i)\le 6\xi $ for $i=1,2$. Applying Lemma \ref{Omega}, we
obtain the inequality $ \dx (s_-, s_+)\le 6L\xi $ for any
component $s$ of $p$. Without loss of generality we may assume
$6L\xi \ge 1$. Thus, $$ \dx (u,v)\le 6L\xi l(p)\le 6L\xi ^2 $$ and
the inequality (\ref{rho}) holds for $\rho =6L\xi ^2$.
\end{proof}


\section{Proof of the main theorem}


We begin by proving the following result, which seems to be of
independent interest. We stress that it does not follow from known
results concerning asymptotic dimensions of hyperbolic graphs
since, in general, the graph $\G $ is not locally finite.

\begin{thm}\label{CG}
The Cayley graph $\G $ has finite asymptotic dimension (with
respect to the metric $\dxh $).
\end{thm}

\begin{proof}
For every $r>0$, we construct a covering of $\G $ as follows. Let
$$A_k=\{ g\in G\; |\; 2kr\le |g|_{X\cup \mathcal H}\le 2(k+1)r\}
$$ and $$S_k=\{ g\in G\; |\; |g|_{X\cup \mathcal H} =2kr\}, $$
where $k=0, 1, \ldots $. To each element $g\in A_k$, $k=1,2,\ldots
$, we associate an arbitrary geodesic $\gamma _g$ in the Cayley
graph $\G $ connecting $g$ to $1$, and denote by $t_g$ the vertex
$\gamma _g\cap S_{k-1}$. Consider the collection $$\mathcal U
(r)=\{ U_k(x)\; | \; k\in \mathbb N, \, x\in S_{k-1} \} ,$$ where
$ U_k(x)=\{ g\in A_k\; |\; t_g=x \} .$

\begin{figure}
  \unitlength 1.2mm 
\linethickness{0.4pt}
\begin{picture}(110,95)(15,0)

\qbezier(19.98,77.78)(65.14,108.19)(109.96,77.78)
\qbezier(27.93,57.98)(66.91,85.65)(102,57.81)
\qbezier(36.95,41.72)(65.85,65.41)(92.98,41.72)
\qbezier(36.24,86.62)(40.84,76.28)(53.92,70.53)
\qbezier(70,71.59)(78.14,74.6)(83.44,90.33)

\put(83.44,76.72){\circle*{1}} \put(75.13,81.49){\circle*{1}}
\put(72.12,52.86){\circle*{1}} \put(81.92,49.27){\circle*{1}}
\put(65.23,2.6){\circle*{1}}
\qbezier(83.44,76.72)(86.53,44.9)(65.23,2.6)
\qbezier(74.95,81.14)(72.74,52.15)(65.23,2.6)

\put(64,20){$\gamma _y$}  \put(75,20){$\gamma $}

\put(63,82.5){\makebox(0,0)[cc]{$U_k(x)$}}
\put(82.2,78.5){\makebox(0,0)[cc]{$a$}}
\put(77.56,82.2){\makebox(0,0)[cc]{$y$}}
\put(89.35,75.51){\makebox(0,0)[cc]{$B(a,r)$}}
\put(38.2,53.96){\makebox(0,0)[cc]{$A_{k-1}$}}
\put(26.91,71.2){\makebox(0,0)[cc]{$A_k$}}
\put(67.78,3.12){\makebox(0,0)[cc]{$1$}}
\put(67,51){\makebox(0,0)[cc]{$x=t_y$}}
\put(80,48.31){\makebox(0,0)[cc]{$s$}}

\put(94.19,76.72){\line(0,1){.575}}
\put(94.18,77.3){\line(0,1){.573}}
\put(94.13,77.87){\line(0,1){.57}}
\put(94.05,78.44){\line(0,1){.565}}
\multiput(93.95,79)(-.0274,.11168){5}{\line(0,1){.11168}}
\multiput(93.81,79.56)(-.03333,.11006){5}{\line(0,1){.11006}}
\multiput(93.64,80.11)(-.03264,.0901){6}{\line(0,1){.0901}}
\multiput(93.45,80.65)(-.03206,.07562){7}{\line(0,1){.07562}}
\multiput(93.22,81.18)(-.03155,.06457){8}{\line(0,1){.06457}}
\multiput(92.97,81.7)(-.03107,.05582){9}{\line(0,1){.05582}}
\multiput(92.69,82.2)(-.03061,.04867){10}{\line(0,1){.04867}}
\multiput(92.39,82.69)(-.03317,.04696){10}{\line(0,1){.04696}}
\multiput(92.05,83.16)(-.03239,.04102){11}{\line(0,1){.04102}}
\multiput(91.7,83.61)(-.03166,.03596){12}{\line(0,1){.03596}}
\multiput(91.32,84.04)(-.03354,.03422){12}{\line(0,1){.03422}}
\multiput(90.92,84.45)(-.03532,.03238){12}{\line(-1,0){.03532}}
\multiput(90.49,84.84)(-.04036,.03321){11}{\line(-1,0){.04036}}
\multiput(90.05,85.21)(-.04208,.03101){11}{\line(-1,0){.04208}}
\multiput(89.58,85.55)(-.04804,.03159){10}{\line(-1,0){.04804}}
\multiput(89.1,85.86)(-.05518,.03219){9}{\line(-1,0){.05518}}
\multiput(88.61,86.15)(-.06393,.03285){8}{\line(-1,0){.06393}}
\multiput(88.1,86.41)(-.07496,.03358){7}{\line(-1,0){.07496}}
\multiput(87.57,86.65)(-.07665,.02953){7}{\line(-1,0){.07665}}
\multiput(87.04,86.86)(-.09113,.02962){6}{\line(-1,0){.09113}}
\multiput(86.49,87.03)(-.11111,.02965){5}{\line(-1,0){.11111}}
\put(85.93,87.18){\line(-1,0){.563}}
\put(85.37,87.3){\line(-1,0){.568}}
\put(84.8,87.39){\line(-1,0){.572}}
\put(84.23,87.45){\line(-1,0){.574}}
\put(83.66,87.47){\line(-1,0){.575}}
\put(83.08,87.47){\line(-1,0){.574}}
\put(82.51,87.43){\line(-1,0){.571}}
\put(81.94,87.37){\line(-1,0){.567}}
\multiput(81.37,87.27)(-.1403,-.0314){4}{\line(-1,0){.1403}}
\multiput(80.81,87.15)(-.11071,-.03111){5}{\line(-1,0){.11071}}
\multiput(80.25,86.99)(-.09074,-.03082){6}{\line(-1,0){.09074}}
\multiput(79.71,86.81)(-.07625,-.03053){7}{\line(-1,0){.07625}}
\multiput(79.18,86.59)(-.0652,-.03024){8}{\line(-1,0){.0652}}
\multiput(78.65,86.35)(-.06349,-.03368){8}{\line(-1,0){.06349}}
\multiput(78.15,86.08)(-.05475,-.03292){9}{\line(-1,0){.05475}}
\multiput(77.65,85.79)(-.04762,-.03221){10}{\line(-1,0){.04762}}
\multiput(77.18,85.46)(-.04167,-.03156){11}{\line(-1,0){.04167}}
\multiput(76.72,85.12)(-.03992,-.03374){11}{\line(-1,0){.03992}}
\multiput(76.28,84.75)(-.03489,-.03284){12}{\line(-1,0){.03489}}
\multiput(75.86,84.35)(-.03308,-.03466){12}{\line(0,-1){.03466}}
\multiput(75.46,83.94)(-.03118,-.03638){12}{\line(0,-1){.03638}}
\multiput(75.09,83.5)(-.03185,-.04144){11}{\line(0,-1){.04144}}
\multiput(74.74,83.04)(-.03255,-.0474){10}{\line(0,-1){.0474}}
\multiput(74.41,82.57)(-.0333,-.05452){9}{\line(0,-1){.05452}}
\multiput(74.11,82.08)(-.03034,-.05622){9}{\line(0,-1){.05622}}
\multiput(73.84,81.57)(-.0307,-.06498){8}{\line(0,-1){.06498}}
\multiput(73.6,81.05)(-.03107,-.07604){7}{\line(0,-1){.07604}}
\multiput(73.38,80.52)(-.03145,-.09052){6}{\line(0,-1){.09052}}
\multiput(73.19,79.98)(-.03188,-.11048){5}{\line(0,-1){.11048}}
\multiput(73.03,79.43)(-.0324,-.14){4}{\line(0,-1){.14}}
\put(72.9,78.87){\line(0,-1){.566}}
\put(72.8,78.3){\line(0,-1){.571}}
\put(72.73,77.73){\line(0,-1){1.149}}
\put(72.69,76.58){\line(0,-1){.575}}
\put(72.71,76.01){\line(0,-1){.572}}
\put(72.76,75.43){\line(0,-1){.569}}
\put(72.85,74.86){\line(0,-1){.563}}
\multiput(72.96,74.3)(.02887,-.11131){5}{\line(0,-1){.11131}}
\multiput(73.1,73.74)(.02898,-.09134){6}{\line(0,-1){.09134}}
\multiput(73.28,73.2)(.02899,-.07685){7}{\line(0,-1){.07685}}
\multiput(73.48,72.66)(.03305,-.07519){7}{\line(0,-1){.07519}}
\multiput(73.71,72.13)(.0324,-.06415){8}{\line(0,-1){.06415}}
\multiput(73.97,71.62)(.0318,-.05541){9}{\line(0,-1){.05541}}
\multiput(74.26,71.12)(.03125,-.04826){10}{\line(0,-1){.04826}}
\multiput(74.57,70.64)(.03071,-.0423){11}{\line(0,-1){.0423}}
\multiput(74.91,70.17)(.03293,-.04059){11}{\line(0,-1){.04059}}
\multiput(75.27,69.73)(.03213,-.03554){12}{\line(0,-1){.03554}}
\multiput(75.66,69.3)(.031369,-.031179){13}{\line(1,0){.031369}}
\multiput(76.06,68.89)(.03574,-.03191){12}{\line(1,0){.03574}}
\multiput(76.49,68.51)(.04079,-.03268){11}{\line(1,0){.04079}}
\multiput(76.94,68.15)(.04673,-.0335){10}{\line(1,0){.04673}}
\multiput(77.41,67.82)(.04845,-.03095){10}{\line(1,0){.04845}}
\multiput(77.89,67.51)(.0556,-.03147){9}{\line(1,0){.0556}}
\multiput(78.39,67.22)(.06435,-.03201){8}{\line(1,0){.06435}}
\multiput(78.91,66.97)(.07539,-.03259){7}{\line(1,0){.07539}}
\multiput(79.44,66.74)(.08987,-.03327){6}{\line(1,0){.08987}}
\multiput(79.98,66.54)(.09152,-.02842){6}{\line(1,0){.09152}}
\multiput(80.52,66.37)(.11149,-.02819){5}{\line(1,0){.11149}}
\put(81.08,66.23){\line(1,0){.564}}
\put(81.65,66.12){\line(1,0){.569}}
\put(82.22,66.04){\line(1,0){.573}}
\put(82.79,65.99){\line(1,0){1.149}}
\put(83.94,65.98){\line(1,0){.573}}
\put(84.51,66.02){\line(1,0){.57}}
\put(85.08,66.09){\line(1,0){.566}}
\multiput(85.65,66.2)(.1398,.0333){4}{\line(1,0){.1398}}
\multiput(86.21,66.33)(.11029,.03256){5}{\line(1,0){.11029}}
\multiput(86.76,66.49)(.09033,.032){6}{\line(1,0){.09033}}
\multiput(87.3,66.68)(.07585,.03153){7}{\line(1,0){.07585}}
\multiput(87.83,66.9)(.0648,.0311){8}{\line(1,0){.0648}}
\multiput(88.35,67.15)(.05604,.03068){9}{\line(1,0){.05604}}
\multiput(88.85,67.43)(.05432,.03363){9}{\line(1,0){.05432}}
\multiput(89.34,67.73)(.0472,.03284){10}{\line(1,0){.0472}}
\multiput(89.81,68.06)(.04125,.0321){11}{\line(1,0){.04125}}
\multiput(90.27,68.41)(.03619,.03141){12}{\line(1,0){.03619}}
\multiput(90.7,68.79)(.03446,.0333){12}{\line(1,0){.03446}}
\multiput(91.12,69.19)(.03263,.03509){12}{\line(0,1){.03509}}
\multiput(91.51,69.61)(.0335,.04013){11}{\line(0,1){.04013}}
\multiput(91.88,70.05)(.0313,.04186){11}{\line(0,1){.04186}}
\multiput(92.22,70.51)(.03192,.04782){10}{\line(0,1){.04782}}
\multiput(92.54,70.99)(.03258,.05495){9}{\line(0,1){.05495}}
\multiput(92.83,71.49)(.0333,.06369){8}{\line(0,1){.06369}}
\multiput(93.1,72)(.02984,.06538){8}{\line(0,1){.06538}}
\multiput(93.34,72.52)(.03007,.07644){7}{\line(0,1){.07644}}
\multiput(93.55,73.05)(.03026,.09092){6}{\line(0,1){.09092}}
\multiput(93.73,73.6)(.03043,.11089){5}{\line(0,1){.11089}}
\multiput(93.88,74.15)(.0306,.1405){4}{\line(0,1){.1405}}
\put(94,74.72){\line(0,1){.568}}
\put(94.1,75.28){\line(0,1){.572}}
\put(94.16,75.85){\line(0,1){.867}}

 \linethickness{0.3pt}

\put(45.04,90.08){\line(1,0){.149}}
\multiput(45.19,90.08)(-.03430391,-.03358924){208}{\line(-1,0){.03430391}}
\multiput(53.07,91.87)(-.034116452,-.033710304){366}{\line(-1,0){.034116452}}
\multiput(60.05,92.76)(-.033727369,-.034039659){476}{\line(0,-1){.034039659}}
\multiput(66,92.91)(-.033954285,-.033686928){556}{\line(-1,0){.033954285}}
\multiput(72.1,92.46)(-.03363128,-.03363128){221}{\line(0,-1){.03363128}}
\multiput(50.99,71.95)(.034259239,.033678574){256}{\line(1,0){.034259239}}
\multiput(77.15,91.87)(-.033719387,-.033719387){626}{\line(0,-1){.033719387}}
\multiput(82.05,90.68)(-.033977201,-.033711754){560}{\line(-1,0){.033977201}}
\multiput(80.57,83.39)(-.034137704,-.033705581){344}{\line(-1,0){.034137704}}
\end{picture}
 \caption{}\label{Ui}
\end{figure}

Obviously the collection $\mathcal U (r) \cup \{ A_0\} $ covers
$G$. Further if $g_1, g_2\in U_k(x)$ for some $k\in \mathbb N$,
$x\in X_k$, we have
$$\dxh (g_1, g_2)\le \dxh (g_1, x)+\dxh (x, g_2)\le 4r. $$
Thus for every $r$, $U(r)$ is uniformly bounded. Finally let
$\sigma $ denote the constant provided by Lemma \ref{xi} and let
$r\ge \sigma $. We consider a ball $B(a,r)=\{ g\in G\; |\; \dxh
(g, a)\le r\} $, where $a\in G$, and assume that $B(a,r)\cap
U_k(x)\ne \emptyset $ for some $k\in \mathbb N$, $ x\in S_{k-1}$.
Let $y\in B(a,r)\cap U_k(x)$. Then $$\dxh (1,a)\ge \dxh (1, y)
-\dxh (y,a)\ge 2kr -r.$$ We fix an arbitrary geodesic $\gamma $ in
$\G $ going from $a$ to $1$ and denote by $s$ the vertex on
$\gamma $ such that $\dxh (1, s)=2kr-2r$. Observe that
$$\dxh (x, y)+\dxh (s, a)\ge 3r
> \dxh (y,a)+\sigma .$$ Therefore, by Lemma \ref{xi} $\dx (x,
s)\le \rho $. Thus for any fixed $k\in \mathbb N$, $B(a, r)$ meets
at most $\mu $ subsets $U _k(x)$,  where $\mu =\sharp\, \{ g\in G,
\; | \; |g|_X\le \rho \} $. Since $B(a,r)$ intersects at most
three annuli $A_k$, we have $\ad \G \le 3\mu $.
\end{proof}

Now we are ready to prove the main result of our paper.

\begin{proof}[Proof of Theorem \ref{main}]
The group $G$ acts on $\G $ by left multiplication. Obviously the
$R$--quasi--stabilizer $W_R(1)$ coincides with the ball $B(R)$ of
radius $R$ with respect to $\dxh $ centered at the identity. It
remains to combine Theorem \ref{CG}, Lemma \ref{Bn}, and Lemma
\ref{ad4}.
\end{proof}


\section{Boundedly generated groups of infinite asymptotic dimension}


Recall that a group $G$ is called {\it boundedly generated} if
there are elements $x_1, \ldots , x_n$ of $G$ such any $g\in G$
can be represented in the form $$g=x_1^{\alpha _1}\ldots
x_n^{\alpha _n}$$ for some $\alpha _1, \ldots , \alpha _n\in
\mathbb Z$. This obviously implies that the Cayley graph $\G $ has
finite diameter for any generating set $X$ and $H_\lambda =\langle
x_\lambda \rangle $, $\lambda \in \{ 1, \ldots , n\} $. Thus
Proposition \ref{prop} is a corollary of the following result.

\begin{prop}
There exists a finitely presented boundedly generated group of
infinite asymptotic dimension.
\end{prop}

\begin{proof}
Recall that a group is called {\it universal} if it contains an
isomorphic copy of any recursively presented group. The existence
of finitely presented universal groups was first proved by Higman
(see \cite[Ch. IV, Theorem 7.3]{LS}). Let $U$ denote a finitely
presented universal group generated by a finite (symmetric) set
$X$.

We set
\begin{equation}\label{fp}
G_1=U\ast \langle a_1\rangle \ast \cdots \ast \langle a_n\rangle .
\end{equation}
Let us enumerate all words $\{ w_1, w_2, \ldots \} $ in the
alphabet $X$ and consider the set
$$
\mathcal R=\{ w_i^{-1} a_1^i\cdots a_n^i,\; i\in \mathbb N\} .
$$
It is easy to see that $\mathcal R$ satisfies the $C^\prime
(\lambda )$ small cancellation condition over the free product
(\ref{fp}), where $\lambda \to 0$ as $n\to \infty $. (For the
definition we refer the reader to \cite[Ch. V, Sec. 9]{LS}.) In
particular, if $n$ is big enough, $\mathcal R$ satisfies $C^\prime
(1/6)$ over  the free product (\ref{fp}) and hence $U$ embeds into
the quotient group
\begin{equation}\label{G2}
G_2=\langle G_1\; |\; R=1,\, R\in \mathcal R\rangle
\end{equation}
(see \cite[Ch. V, Corollary 9.4]{LS}).

Note that the presentation (\ref{G2}) is recursive. As $U$ is
universal, it contains an isomorphic copy of $G_2$. Let $\alpha
\colon G_2\to U$ be the monomorphism that maps $G_2$ to its copy
in $U\le G_2$. By $G$ we denote the ascending HNN--extension
\begin{equation}\label{G3}
G=\langle G_2,t\; |\; g^t=\alpha (g), \, g\in G_2 \rangle .
\end{equation}

Let $g$ be an arbitrary element of $G$. Observe that the subgroup
$$
N=\bigcup\limits_{j=1}^\infty t^jG_2t^{-j}=
\bigcup\limits_{j=1}^\infty t^jUt^{-j}
$$
is normal in $G$ and thus the kernel of the natural homomorphism
$G\to \langle t\rangle $ coincides with $N$. Therefore there exist
$j,k\in \mathbb N$ such that $t^{-j}(gt^k)t^j\in U$. Furthermore,
any element of $U$ can be represented as a product $a_1^i\cdots
a_n^i$ for a certain $i$ according to the relations $R=1$, $R\in
\mathcal R$. Hence any element $g\in G$ can be represented as
$$g=t^{\alpha _1} a_1^i\cdots a_n^it^{\alpha _2} $$ for some
$\alpha _1$, $\alpha _2$, $i \in \mathbb Z$. In particular, the
group $G$ is boundedly generated.

One can also observe that $G$ is finitely presented. Indeed
expanding the presentation (\ref{G3}) we obtain
\begin{equation}\label{G4}
G=\langle a_1, \ldots a_n,t, U\; |\; R=1,\, R\in \mathcal R, \;
g^t=\alpha (g), \, g\in G_2 \rangle .
\end{equation}
Since $(G_2)^t\le U$ and words from the set $\mathcal R$ involve
elements of $G$ only, all relations of the form $R=1$, $R\in
\mathcal R$, follow from the relations of the group $U$. Thus we
can omit the set of relations $R=1$, $R\in \mathcal R$ in
(\ref{G4}) and get a finite presentation for $G$. Finally we
notice that $\ad G=\infty $ since any universal group contains an
isomorphic copy of $\mathbb Z^m$ for all $m\in \mathbb N$ and $\ad
\mathbb Z^m=m$.
\end{proof}


\begin{thebibliography}{99}

\bibitem{BD1}
G. Bell, A. Dranishnikov, {\it On asymptotic dimension of groups},
Alg. and Geom. Topology, {\bf 1} (2001), 57--71.

\bibitem{BD}
G. Bell, A. Dranishnikov, {\it On asymptotic dimension of groups
acting on trees}, Geom. Dedicata {\bf 103} (2004), 89--101.

\bibitem{BD2}
G. Bell, A. Dranishnikov, {\it A Hurewicz-type theorem for
asymptotic dimension and applications to geometric group theory},
prep., 2004; available at
http://www.arxiv.org/abs/math.GR/0407431.

\bibitem{Bow}
B.H. Bowditch, {\it Relatively hyperbolic groups,} prep., 1999.

\bibitem{BriH}
M. Bridson, A. Haefliger, Metric spaces of non--positive
curvature, Springer, 1999.

\bibitem{CG}
G. Carlsson, B. Goldfarb, {\it The integral $K$-theoretic Novikov
conjecture for groups with finite asymptotic dimension}, Invent.
Math. {\bf 157} (2004), no. 2, 405--418.

\bibitem{Dah}
F. Dahmani, {\it Combination of convergence groups,} Geom. Topol.
{\bf 7} (2003), 933--963 (electronic).

\bibitem{DY}
F. Dahmani, A. Yaman, {\it Bounded geometry in relatively
hyperbolic groups}, prep., 2004.

\bibitem{D}
A. Dranishnikov, {\it On hypersphericity of manifolds with finite
asymptotic dimension}, Trans. Amer. Math. Soc. {\bf 355} (2003),
no. 1, 155--167.

\bibitem{SD} C. Drutu, M.V. Sapir, {\it Tree graded spaces
and asymptotic cones} prep., 2004; available at
http://www.arxiv.org/abs/math.GT/0405030.

\bibitem{F}
B. Farb, {\it Relatively hyperbolic groups,} GAFA, {\bf 8} (1998),
810--840.

\bibitem{GhH}
E. Ghys, P. de la Harpe, Eds., Sur les groupes hyperboliques
d'apr\'es Mikhael Gromov, Progress in Math., 83, Birka\"user,
1990.

\bibitem{Gro}
M. Gromov, {\it Hyperbolic groups,} Essays in Group Theory, MSRI
Series, Vol.8, (S.M. Gersten, ed.), Springer, 1987, 75--263.

\bibitem{LS}
R.C. Lyndon, P.E. Shupp, Combinatorial Group Theory,
Springer--Verlag, 1977.

\bibitem{WH} D.V. Osin, {\it Weak hyperbolicity and free constructions,}
Contemp. Math., {\bf 360} (2004), 103--111.

\bibitem{RHG} D.V. Osin, {\it Relatively hyperbolic groups: Intrinsic
geometry, algebraic properties, and algorithmic problems}, Memoirs
of Amer. Math. Soc., to appear.

\bibitem{R}
J. Roe, Lectures on coarse geometry. University Lecture Series,
31. American Mathematical Society, Providence, RI, 2003

\bibitem{Yaman}
A. Yaman, {\it A topological characterization of relatively
hyperbolic groups}, J. Reine Angew. Math. {\bf 566} (2004),
41--89.

\bibitem{Yu}
G. Yu, {\it The Novikov conjecture for groups with finite
asymptotic dimension}, Ann. of Math. (2) {\bf 147} (1998), no. 2,
325--355.

\bibitem{Yu1}
G. Yu, {\it Zero-in-the-spectrum conjecture, positive scalar
curvature and asymptotic dimension}, Invent. Math. {\bf 127}
(1997), no. 1, 99--126.

\end{thebibliography}
\end{document}